%% file: sat.tex
\documentclass{amsart}
\usepackage{graphicx}
\usepackage{psfrag}

\input{preamble}

\title[The saturation conjecture (after A.~Knutson and T.~Tao)]{The
saturation conjecture \\ (after A.~Knutson and T.~Tao)}
\date{\today}
\author{Anders Skovsted Buch}
\address{Department of Mathematics\\
  University of Chicago\\
  Chicago, IL 60637}
\email{abuch@math.uchicago.edu} 

\input{psfrag}

\begin{document}

\maketitle

The purpose of this exposition\footnote{given as a talk at UC
Berkeley, September 1998} is to give a simple and complete treatment
of Knutson and Tao's recent proof of the saturation conjecture
\cite{knutson.tao:honeycomb}.

If $\lambda$ is a partition of length at most $n$, let $V_\lambda$
denote the corresponding highest weight representation of $\GL_n(\C)$.
Define
\[ T_n = \{ (\lambda,\mu,\nu) \mid V_\nu \subset V_\lambda \otimes
V_\mu \} \,.
\]
This set is important in numerous areas besides representation theory.
In Schubert calculus it describes when an intersection of Schubert
cells must be non-empty.  In combinatorics, a triple is in $T_n$ if
and only if there exists a Littlewood-Richardson skew tableau with
shape $\nu/\lambda$ and content $\mu$.

It is well known that $T_n \subset \Z^{3n}$ is a semi-group under
addition, a fact which Zelevinsky attributes to Brion and Knop
\cite{zelevinsky:littlewood-richardson}.  Klyachko has given
\cite{klyachko:stable} a nice description of the saturation
\[ \bar T_n = \{ (\lambda,\mu,\nu) \mid \exists N \in \N : 
   (N\lambda, N\mu, N\nu) \in T_n \} \,.
\]
A triple $(\lambda,\mu,\nu)$ is in $\bar T_n$ if and only if the
entries of $\lambda$, $\mu$, and $\nu$ satisfy certain inequalities
that come from Schubert calculus (see also \cite{fulton:eigenvalues}).
This made the following conjecture particularly important.

\begin{satconj}
$(\lambda,\mu,\nu) \in T_n \Longleftrightarrow
(N\lambda, N\mu, N\nu) \in T_n$.
\end{satconj}
In other words $T_n$ is saturated in $\Z^{3n}$.  Note that the
implication $\Rightarrow$ is a trivial consequence of the fact
that $T_n$ is a semi-group or of the original Littlewood-Richardson
rule.

In July 1998, Knutson and Tao gave a proof of this conjecture, using
two wonderful new descriptions of Berenstein-Zelevinsky polytopes
called the {\em honeycomb\/} and {\em hive\/} models
\cite{knutson.tao:honeycomb}.  The goal of this exposition is to
present a simple and complete proof based only on the hive model.
Since the final version of Knutson and Tao's paper will likely be
based on honeycombs alone, we hope that this may be useful.  Most
constructions used here come directly from the first version of
Knutson and Tao's preprint, even if they may be replaced by honeycomb
equivalents in their published paper.  One innovation, in Section
\ref{sec_smallflat}, is the construction of a graph from a hive, which
is used to simplify their argument.  In an appendix of Fulton it is
shown that the hive model is equivalent to the original
Littlewood-Richardson rule.

\section{The hive model}

We start with a triangular array of {\em hive vertices}, $n+1$ on
each side.
\[ \picA{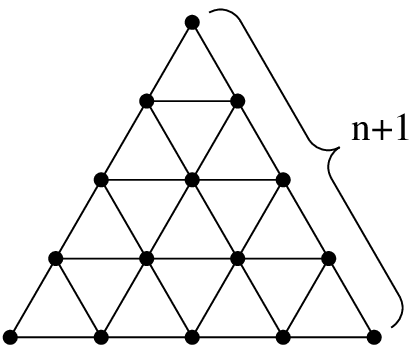} \]
This array is called the (big) {\em hive triangle\/}.  When lines are
drawn through the hive vertices as shown, the hive triangle is split
up into $n^2$ {\em small triangles}.  By a {\em rhombus\/} we mean the
union of two small triangles next to each other.

Let $H$ be the set of hive vertices and $\R^H$ the labelings of these
by real numbers.  Each rhombus gives rise to an inequality on $\R^H$
saying that the sum of the labels at the obtuse vertices must be
greater than or equal to the sum of the labels at the acute vertices.
\[ \picA{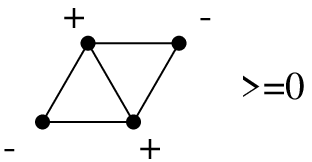} \]
A {\em hive\/} is a labeling that satisfies all rhombus inequalities.
A hive is {\em integral\/} if all its labels are integers.  We let $C
\subset \R^H$ denote the convex polyhedral cone consisting of all
hives.

The Littlewood-Richardson coefficient $c^\nu_{\lambda\mu}$ is defined
to be the multiplicity of $V_\nu$ in $V_\lambda \otimes V_\mu$.
Equivalently it is the number of Littlewood-Richardson skew tableaux
of shape $\nu/\lambda$ and content $\mu$ \cite[\S 5.2]{fulton:young}.
The following theorem illustrates the importance of
Berenstein-Zelevinsky polytopes.
\begin{thm}
\label{thm_lrrule}
Let $\lambda$, $\mu$, and $\nu$ be partitions with $|\nu| = |\lambda|
+ |\mu|$.  Then $c^\nu_{\lambda\mu}$ is the number of integral hives
with border labels:
\[ \picA{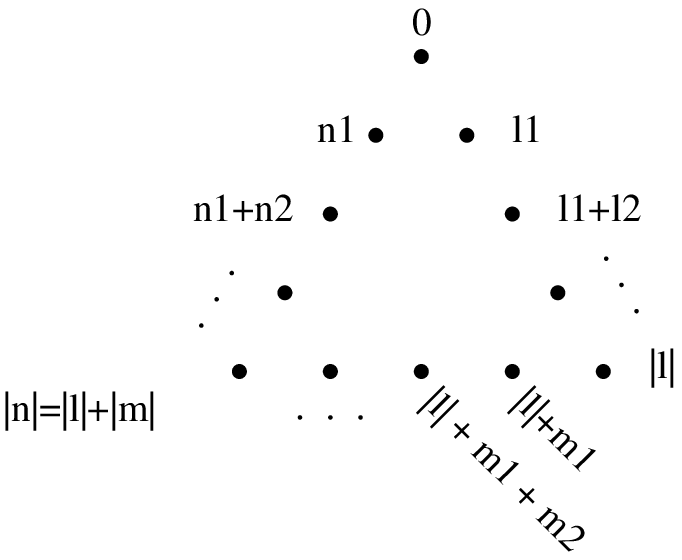} \]
\end{thm}
Knutson and Tao prove this by translating hives with integer labels
into tail-positive Berenstein-Zelevinsky patterns, which are known to
count $c^\nu_{\lambda\mu}$ \cite{berenstein.zelevinsky:triple},
\cite{zelevinsky:littlewood-richardson}.  An alternative direct proof
of Fulton can be found in the appendix.

\begin{example}
To compute $c^{3\,2\,1}_{2\,1\,,\,2\,1}$ we can take $n=3$ and border
labels as in the picture.
\[ \picA{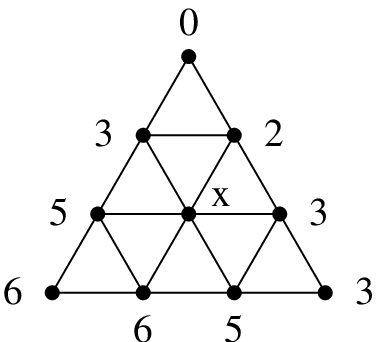} \]
Let $x$ be the undetermined label of the middle hive vertex.  Then the
rhombus inequalities say that $4 \leq x \leq 5$.  It follows that
there are two integral hives with this border, so
$c^{3\,2\,1}_{2\,1\,,\,2\,1} = 2$.
\end{example}

Let $B$ be the set of border vertices, and $\rho : \R^H \to \R^B$ the
restriction map.  The restriction of a hive to the border vertices by
$\rho$ is called its {\em border}.  For $b \in \R^B$, the fiber
$\rho^{-1}(b) \cap C$ is easily seen to be a compact polytope, which
we will call the {\em hive polytope\/} over $b$.  If $b$ comes from a
triple of partitions as in \refthm{thm_lrrule}, this is also called
the hive polytope over the triple.  We will call the vertices of a
hive polytope for its {\em corners}.

We can now describe the strategy of Knutson and Tao's proof.  If
$(N\lambda, N\mu, N\nu)$ is in $T_n$, then the hive polytope over this
triple contains an integral hive.  By scaling this polytope down by a
factor $N$, it follows that the hive polytope over $(\lambda,\mu,\nu)$
is not empty.  Therefore it is enough to show that if $b \in \Z^B$ and
$\rhoinv(b) \cap C \neq \emptyset$ then $\rhoinv(b) \cap C$ contains
an integral hive.

Let $\omega$ be a functional on $\R^{H-B}$ which maps a hive to a
linear combination of the labels at non-border vertices, with generic
positive coefficients.  Then for each $b \in \rho(C)$, $\omega$ takes
its maximum at a unique hive in $\rhoinv(b) \cap C$.  The strategy is
to prove that this hive is integral if $b$ is integral.

\begin{example}
\label{example_nonint}
Even though all rhombus inequalities are integrally defined, a hive
polytope over an integral border can still have non-integral corners.
The following hive is an example, and therefore it does not maximize
any generic positive functional $\omega$.
\[ \picA{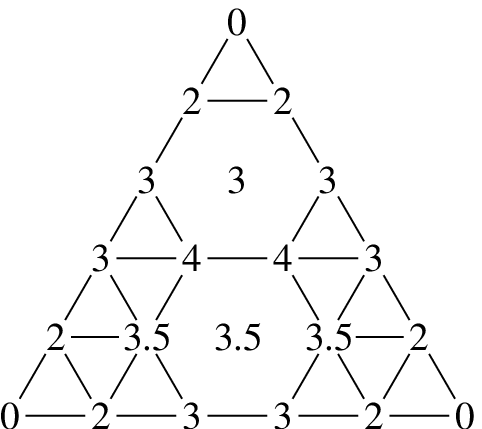} \]
In the picture we have omitted the lines across rhombi where the
rhombus inequality is satisfied with equality, which makes it easy to
see that this hive is a corner of its hive polytope.  It is not hard
to show that for $n \leq 4$ and $b \in \Z^B$, all corners of
$\rhoinv(b) \cap C$ are integral hives.
\end{example}

\section{Flatspaces}

We can consider a hive as a graph over the hive triangle.  At each
hive vertex we use the label as the height.  We then extend these
heights to a graph over the entire hive triangle by using linear
interpolation over each small triangle.  A rhombus inequality now says
that the graph over the rhombus cannot bend up across the middle line.
\[ \picA{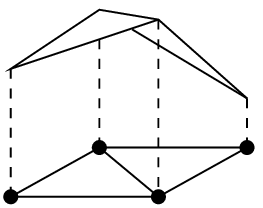} \]
In this way the graph becomes the surface of a convex mountain.  The
graph is flat (but not necessarily horizontal) over a rhombus if and
only if the rhombus inequality is satisfied with equality.

We define a {\em flatspace\/} to be a maximal connected union of small
triangles such that any contained rhombus is satisfied with equality.
The flatspaces split the hive triangle up in disjoint regions over
which the mountain is flat.  The flatspaces of the hive in
\refexample{example_nonint} consist of two hexagons and 13 small
triangles.

Flatspaces have a number of nice properties.  We will list the ones we
need below.  Since all of these are straightforward to prove directly
from the definitions, we will simply give intuitive reasons for them.

\noindent {\bf 1.  Flatspaces are convex.}  This is clear since they
lie under intersections of a convex mountain with a (convex) plane.

\noindent {\bf 2.  All flatspaces have one of the following five
shapes (up to rotations and different side lengths):}
\[
\picB{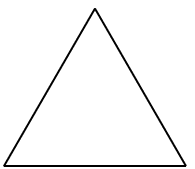} \hspace{0.5cm}
\picB{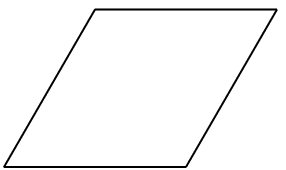} \hspace{0.5cm}
\picB{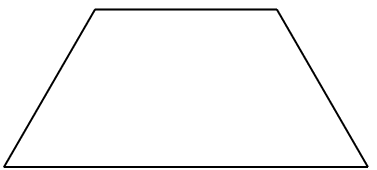} \hspace{0.5cm}
\picB{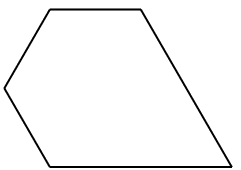} \hspace{0.5cm}
\picB{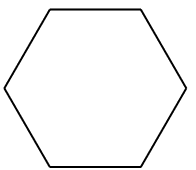}
\]
These are the only convex shapes that can be constructed from small
triangles.

\noindent {\bf 3.  A side of a flatspace is either on the border of
the big hive triangle, or it is also a side of a neighbor flatspace.}
In other words, a side of one flatspace can't be shared between
several neighbor flatspaces.  This again follows from the convexity of
the mountain described above.

Given a labeling $b \in \R^B$, let $x, y, z$ be labels of
consecutive border vertices on the same side of the big hive triangle.
\[ \picA{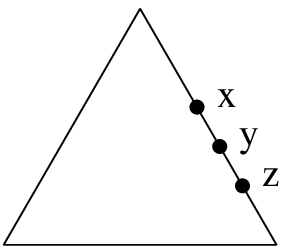} \]
If $b$ is the border of a hive, then the rhombus inequalities imply
that $y-x \geq z-y$.  We will say that $b$ is {\em regular\/} if we
always have $y-x > z-y$, when $x,y,z$ are chosen in this way.  Regular
borders correspond to triples of strictly decreasing partitions.

\noindent {\bf 4.  If the border of a hive is regular then no
flatspace has a side of length $\geq 2$ on the border of the big hive
triangle.}  In fact, if the labels $x,y,z$ above are on a flatspace
side, then $y-x = z-y$.

Given a hive, a non-empty subset $S \subset H-B$ is called
increasable, if the same small positive amount can be added to the
labels of all hive vertices in $S$, such that the labeling is still a
hive.

\noindent {\bf 5.  The interior vertices of a hexagon-shaped flatspace
is an increasable subset.}  Proving this is a matter of checking that
each rhombus inequality still holds after adding a small enough amount
to the labels of these vertices.  Only rhombi that are already flat
need to be considered, since for all others there is some ``slack to
cut''.

Note that the corresponding statements for flatspaces of other shapes
are wrong.  The reason is that all other shapes have at least one
sharp corner (with a $60^\circ$ angle).  Lifting the interior vertex
closest to a sharp corner is prohibited by the inequality of the
rhombus in that corner.
\[ \picA{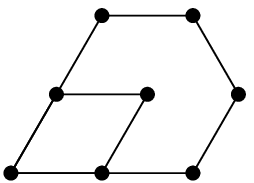} \]

\begin{prop}
\label{prop_smallflat}
If a hive with regular border has no increasable subsets, then its
flat\-spaces consist of small triangles and small rhombi.
\end{prop}
\begin{proof}
Otherwise some flatspace has a side of length $\geq 2$.  This follows
because the only types of flatspaces that have all sides of length one
are small triangles, rhombi, and small hexagons, and the later do not
occur by property 5.

Let $m$ be the maximal length among all sides of flatspaces.  We will
proceed by constructing a region consisting of flatspaces with a side
of length $m$, such that the interior hive vertices of the region is
an increasable subset.  The crucial point is to avoid sharp corners
pointing out from the region.  We need $m \geq 2$ to be sure that
interior vertices exist.

Start by taking any flatspace having a side of length $m$, and mark
this side.  In the pictures this is shown by making the side thick.
Then choose a line crossing (the extension of) the marked side in an
angle of $60^\circ$ and call it the moving direction.  If the
flatspace is a triangle or a parallelogram, we furthermore mark an
additional side.  For a triangle, this is the other side not parallel
to the moving direction, while for a parallelogram we mark the side
opposite the one already marked.
\[ \picB{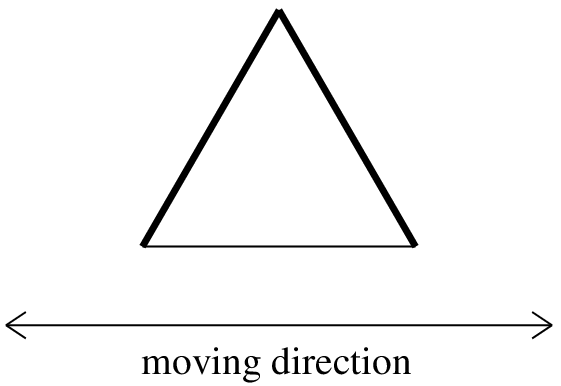} \]

We construct a region, starting with the chosen flatspace.  As long as
the region has a marked side on its outer border, the flatspace on the
opposite side is added to the region.  Note that this flatspace is
well defined by property 3, since regularity prevents any marked edges
from being on the border of the big hive triangle.  If the new
flatspace is a triangle, we mark its unmarked side which is not
parallel to the moving direction.  If the new flatspace is a
parallelogram, we mark the side opposite the old marked side.
\[ \picB{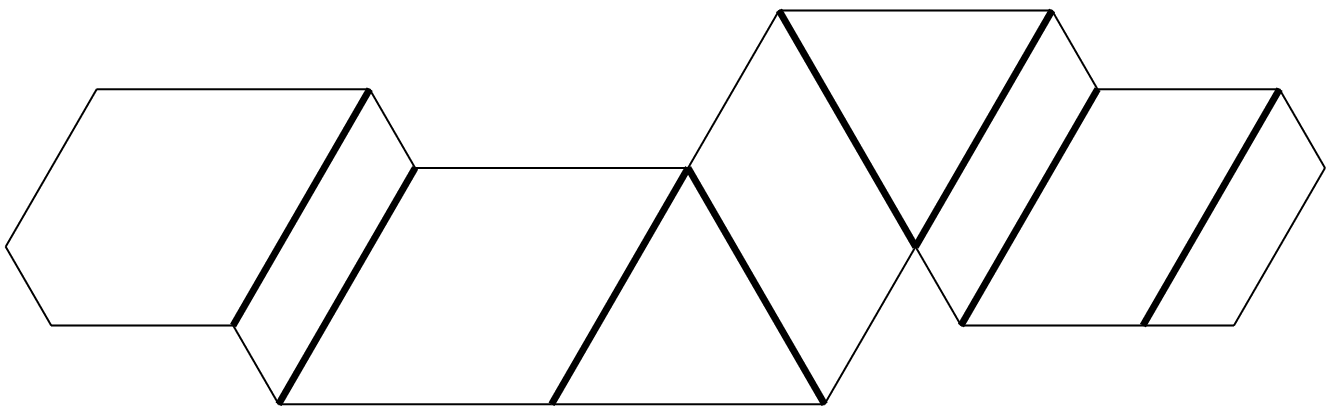} \]

Since the region always grows along the moving direction, it will
never go in loops.  Now since no marked edge can ever reach the border
of the big hive triangle, the described process will stop.  By
induction it is easy to see that all sharp corners must be on marked
sides, and since the final region has no marked sides on its boundary,
this region has no outward sharp corners.  It is now easy to verify
that the interior vertices form an increasable subset, and this
contradiction completes the proof.
\end{proof}

\section{Small Flatspaces}
\label{sec_smallflat}

Let $h$ be a hive, all of whose flatspaces are small triangles or
small rhombi.  We construct a (colored) graph $G$ from $h$ as
follows.  $G$ has one blue vertex in the middle of each small triangle
flatspace.  In addition there is one red vertex on every flatspace
side.  Each blue vertex is connected to the three vertices on the
sides of its triangle, and the two red vertices on opposite sides of
a rhombus are connected.  This graph is topologically equivalent to
the reduced honeycomb tinkertoy of Knutson and Tao.
\[ \picA{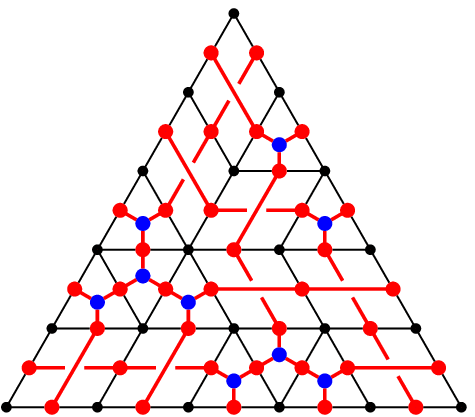} \]

\begin{lemma}
\label{lemma_acyclic}
If $h$ is a corner of its hive polytope $\rho^{-1}(\rho(h)) \cap C$,
then $G$ is acyclic.
\end{lemma}
\begin{proof}
Suppose $G$ has a non-trivial loop, and give this loop an orientation.
Each hive vertex then has a well defined winding number, which is the
number of times the loop goes around this vertex, counted positive in
the counter clockwise direction.  Note that the winding number is zero
for each border vertex, and that some winding numbers are non-zero if
the loop is not trivial.

For each $r \in \R$, let $h_r \in \R^H$ be the labeling which maps
each hive vertex to the label of $h$ at the vertex plus $r$ times the
winding number of the vertex.  We will show that $h_r$ is a hive for
$r \in (-\epsilon, \epsilon)$, $\epsilon > 0$.  This implies that $h$
is an interior point of a line segment contained in its hive
polytope, which contradicts the assumption that $h$ is a corner.

Choose any $\epsilon > 0$ such that each rhombus inequality that is
strict for $h$ is also satisfied for $h_r$ when $|r| < \epsilon$.  We
claim that this $\epsilon$ will do.  Consider any rhombus satisfied by
$h$ with equality: 
\[ \picA{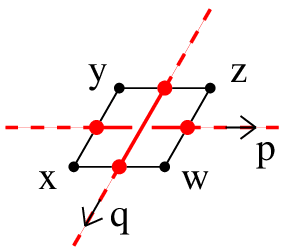} \]
Suppose that the loop goes through the two edges with multiplicities
$p$ and $q$ in the indicated directions, and that the vertex with
label $x$ has winding number $t$.  Then going clockwise around the
rhombus, the winding numbers of the three other vertices are $t + p$,
$t + p + q$, and $t + q$.  It follows that the labels of $h_r$ are
\begin{align*}
y' &= y + r (t + p) & z' &= z + r (t + p + q) \\
x' &= x + r t & w' &= w + r (t + q)
\end{align*}
Since the rhombus is flat for $h$, we have $x + z = y + w$.  But this
implies that $x' + z' = y' + w'$, and so the rhombus is also flat for
$h_r$.
\end{proof}

\begin{prop}
\label{prop_intlincomb}
Let $h$ be a hive which is a corner of its hive polytope
$\rho^{-1}(\rho(h)) \cap C$.  Suppose the flatspaces of $h$ consist
only of small triangles and small rhombi.  Then the labels of $h$ are
integer linear combinations of the border labels.
\end{prop}
\begin{proof}
By \reflemma{lemma_acyclic}, the graph $G$ for $h$ is acyclic.  Label
each red vertex with the difference of the labels of the hive
vertices on its side:
\[ \picA{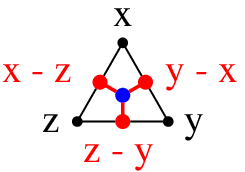} \]

By construction, the sum of the labels of three red vertices
surrounding any blue vertex is zero.  Furthermore, if two red vertices
are connected by a single edge, then their labels are equal.  This
follows because the rhombus that separates them is satisfied with
equality.  We will show that all red vertex labels are $\Z$-linear
combinations of the border labels.  Since this implies that also all
labels of hive vertices are such linear combinations, this will finish
the proof.

Consider any connected component of $G$.  We will use the induction
hypothesis that each red endpoint vertex label is a $\Z$-linear
combination of the border labels.  This is clearly true at the
starting point, since all endpoint vertices are on the border of the
big hive triangle.

If an endpoint vertex is connected to another red vertex, we can
remove the endpoint (and its edge) from the graph, making the other
red vertex a new endpoint.  Since the two vertex labels are equal, the
induction hypothesis remains valid.

If no endpoint vertices are connected to a red vertex, then since $G$
is acyclic, some blue vertex must be connected to two endpoint
vertices.  Since the third red vertex connected to this blue vertex
has a label which is minus the sum of the two other, it is a
$\Z$-linear combination of the border labels.  Therefore we can
discard the blue vertex and the two endpoint vertices (with all their
edges), without violating the induction hypothesis.  Continuing in
this way, the whole graph is eventually removed, and at this point we
have verified that all labels are $\Z$-linear combinations of the
border labels.
\end{proof}

\section{Proof of the saturation conjecture}

We will call a functional on $\R^{H-B}$ is {\em generic\/} if it takes
its maximum at a unique point in $\rhoinv(b) \cap C$ for each $b \in
\rho(C)$.  The existence of generic functionals follows from the
existence of secondary fans in linear programming \cite[\S
1]{sturmfels.thomas:variation}.  We can now finish the proof of the
saturation conjecture.

\begin{thm}
$(\lambda,\mu,\nu) \in T_n \Longleftrightarrow 
(N\lambda, N\mu, N\nu) \in T_n$.
\end{thm}
\begin{proof}
As already noted, it is enough to show that if $b \in \rho(C) \cap
\Z^B$ then the fiber $\rhoinv(b) \cap C$ contains an integral hive.

Fix a generic functional $\omega$ on $\R^{H-B}$ which maps a hive to a
linear combination with positive coefficients of the labels at
non-border hive vertices.  For each $b \in \rho(C)$, let $\ell(b)$ be
the unique point in $\rhoinv(b) \cap C$ where $\omega$ is maximal.
Then $\ell : \rho(C) \to C$ is a continuous piece-wise linear map
\cite[\S 1]{sturmfels.thomas:variation}.  Notice that $\ell(b)$ has no
increasable subsets.

For a regular border $b \in \rho(C)$, \refprop{prop_smallflat} implies
that the flatspaces of $\ell(b)$ consist of small triangles and
rhombi.  By \refprop{prop_intlincomb} this implies that all labels of
$\ell(b)$ are $\Z$-linear combinations of the labels of $b$.  By
continuity it follows that each linear piece of $\ell$ is integrally
defined.  In particular $\ell(b)$ is an integral hive if $b$ is
integral.
\end{proof}

\section{Remarks and questions}

Knutson and Tao's proof of the saturation conjecture implies that
Klyachko's inequalities for $T_n$ can be produced by a simple
recursive algorithm, which uses the inequalities for $T_k$, $1 \leq k
\leq n-1$ \cite{klyachko:stable}, \cite{knutson.tao:honeycomb},
\cite{fulton:eigenvalues}.  Another important consequence is 
Horn's conjecture, which says that the same inequalities describe
which sets of eigenvalues can arise from two Hermitian matrices and
their sum \cite{horn:eigenvalues}.

In connection with Klyachko's work, it has been of interest which
triples $(\lambda,\mu,\nu)$ have Littlewood-Richardson coefficient
$c^\nu_{\lambda\mu}$ equal to one.  Fulton has conjectured that this
is equivalent to $c^{N\nu}_{N\lambda, N\mu}$ being one for any $N \in
\N$.  This has been verified in all cases with $N |\nu| \leq 66$.

For $n = 3$ it is easy to show that a triple of partitions has
coefficient one if and only if it corresponds to a point on the
boundary of the cone $\rho(C)$.  In general, Fulton's conjecture
implies that the triples with coefficient one are exactly those
corresponding to points in a collection of faces of $\rho(C)$.  For $n
\geq 3$ this means that all triples corresponding to interior points
in $\rho(C)$ have coefficient at least two.

One approach for proving Fulton's conjecture is to show that if $b \in
\rho(C) \cap \Z^B$, then any generic positive functional $\omega$ on
$\R^{H-B}$ must be minimized (as well as maximized) at an integral
hive in $\rhoinv(b) \cap C$.  In fact, by \refprop{prop_intlincomb} it
is enough to prove:

{\em If $b \in \rho(C)$ is a generic border and if a generic positive
functional $\omega$ is minimized at $h \in \rhoinv(b) \cap C$, then
the flatspaces of $h$ consist of small triangles and rhombi.}

Part of proving this is to specify when a border $b$ is generic.  We
believe the statement is true if $b$ avoids finitely many hyperplanes
in $\R^B$.

The Littlewood-Richardson coefficients $c^\nu_{\lambda\mu}$ have the
following natural generalization.  Given decreasing sequences of
integers $\nu$, and $\lambda(1), \dots, \lambda(r)$, let
$c^\nu_{\lambda(1),\dots,\lambda(r)}$ denote the multiplicity of
$V_\nu$ in the holomorphic representation $V_{\lambda(1)} \otimes
\cdots \otimes V_{\lambda(r)}$.  When $\nu = (0,\dots,0)$, this
specializes to the symmetric Littlewood-Richardson coefficient
$c_{\lambda(1),\dots,\lambda(r)}$ which is the dimension of the
$\GL_n(\C)$-invariant subspace of $V_{\lambda(1)} \otimes \cdots
\otimes V_{\lambda(r)}$.  Postnikov and Zelevinsky have pointed out
that the saturation conjecture as stated in the introduction implies a
similar result for these generalized coefficients, i.e.\ 
\comment{same line?}
\begin{equation}
\label{eqn_nsat}
c^\nu_{\lambda(1),\dots,\lambda(r)} \neq 0 \Longleftrightarrow
c^{N\nu}_{N\lambda(1), \dots, N\lambda(r)} \neq 0 \,.
\end{equation}
Knutson has shown us that by combining several hive triangles, one
obtains a polytope whose integral points count these more general
coefficients.  This gives rise to another proof of \refeqn{eqn_nsat}.

In \cite{buch.fulton:chern} another type of generalized
Littlewood-Richardson coefficients 
related to quiver varieties are described.  It would be very
interesting if these coefficients can be realized as the number of
integral points in some polytope.

\appendix

\input{appendix}

\input{bibliography}

\end{document}

%% file: preamble.tex
\newcommand{\N}{{\mathbb N}}
\newcommand{\Z}{{\mathbb Z}}

\newcommand{\R}{{\mathbb R}}
\newcommand{\C}{{\mathbb C}}

\newcommand{\GL}{\operatorname{GL}}
\newcommand{\rhoinv}{\rho^{-1}}
\newcommand{\tth}{{\text{th}}}

\newcommand{\picA}[1]{\includegraphics{#1}}
\newcommand{\picB}[1]{\includegraphics[scale=0.75]{#1}}

\theoremstyle{plain}
\newtheorem*{satconj}{Saturation conjecture}
\newtheorem{thm}{Theorem}
\newtheorem{lemma}{Lemma}
\newtheorem{prop}{Proposition}

\theoremstyle{definition}
\newtheorem{example}{Example}

\numberwithin{equation}{section}

\newcommand{\refeqn}[1]{(\ref{#1})}
\newcommand{\refthm}[1]{Theorem~\ref{#1}}
\newcommand{\refprop}[1]{Proposition~\ref{#1}}
\newcommand{\reflemma}[1]{Lemma~\ref{#1}}

\newcommand{\refexample}[1]{Example~\ref{#1}}

\newcommand{\comment}[1]{}

%% file: psfrag.tex
\psfrag{n+1}{$n+1$}
\psfrag{>=0}{$\geq 0$}
\psfrag{n1}{$\nu_1$}
\psfrag{n1+n2}{$\nu_1 + \nu_2$}
\psfrag{|n|=|l|+|m|}{$|\nu| = |\lambda| + |\mu|$}
\psfrag{l1}{$\lambda_1$}
\psfrag{l1+l2}{$\lambda_1 + \lambda_2$}
\psfrag{|l|}{$|\lambda|$}
\psfrag{|l|+m1}{$|\lambda| + \mu_1$}
\psfrag{|l| + m1 + m2}{$|\lambda| + \mu_1 + \mu_2$}
\psfrag{x}{$x$}
\psfrag{y}{$y$}
\psfrag{z}{$z$}
\psfrag{w}{$w$}
\psfrag{p}{$p$}
\psfrag{q}{$q$}
\psfrag{y - x}{$y - x$}
\psfrag{z - y}{$z - y$}
\psfrag{x - z}{$x - z$}
\psfrag{a10}{$a^1_0$}
\psfrag{a11}{$a^1_1$}
\psfrag{a12}{$a^1_2$}
\psfrag{a13}{$a^1_3$}
\psfrag{a20}{$a^2_0$}
\psfrag{a21}{$a^2_1$}
\psfrag{a22}{$a^2_2$}
\psfrag{a30}{$a^3_0$}
\psfrag{a31}{$a^3_1$}
\psfrag{a40}{$a^4_0$}
\psfrag{ai+1k}{$a^{i+1}_k$}
\psfrag{aik+1}{$a^i_{k+1}$}
\psfrag{aik-1}{$a^i_{k-1}$}
\psfrag{ai+1k-1}{$a^{i+1}_{k-1}$}
\psfrag{aik}{$a^i_k$}
\psfrag{ai-1k}{$a^{i-1}_k$}

%% file: appendix.tex
\section{A direct proof of \refthm{thm_lrrule}}

The aim of this appendix is to give a simple and direct bijection
between the hives with given boundary (given by a triple of
partitions), and the set of Littlewood-Richardson skew tableaux for
the given triple.  In principle one could construct such a mapping
from \cite{carre:rule}, but it is simpler to do it directly from
hives; in the description we give here, it is easy to see that the map
is a bijection, without knowing that the two sets have the same
cardinality. As in \cite{carre:rule}, we produce contratableaux, but
there is a standard bijection between these and the original
Littlewood-Richardson skew tableaux.

Consider an integral hive, with sides having $n+1$ entries,
corresponding to partitions $\lambda$, $\mu$, and $\nu$, with $|\nu| =
|\lambda| + |\mu|$.  The differences down the northwest to southeast
border give the partition $\lambda$, the differences across the bottom
border from right to left give $\mu$, and the differences down the
northeast to southwest border give $\nu$ (see \refthm{thm_lrrule}).
The main idea for constructing a skew tableau with a reverse-lattice
word is to use the other northwest to southeast rows of entries to
construct a chain of subpartitions of $\lambda$.

The entries of the hive will be denoted $a^i_k$, with $1 \leq i \leq
n+1$ and $0 \leq k \leq n+1-i$. Here the superscript denotes the
northwest to southeast row of the entry, with the first row being the
long row on the boundary, and the others in order below that; the
subscripts number the entries along the rows, from northwest to
southeast.
\[ \picA{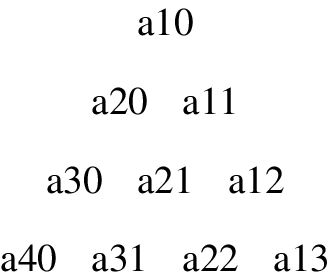} \]
Note that $a^1_0 = 0$, and that $\lambda_k = a^1_k - a^1_{k-1}$ for $1
\leq k \leq n$.

For $1 \leq i \leq n$ define a sequence $\lambda^{(i)} =
(\lambda^{(i)}_1, \dots, \lambda^{(i)}_{n+1-i})$ by setting
$\lambda^{(i)}_k = a^i_k - a^i_{k-1}$.  Note that $\lambda^{(1)} =
\lambda$.

There are three types of rhombus inequalities, depending on the
orientation of the rhombus.  We first consider two of them:

\vspace{0.2cm} \noindent
(1) \raisebox{-0.3cm}{\picB{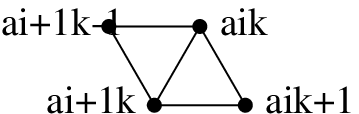}} This says that
$\lambda^{(i+1)}_k \geq \lambda^{(i)}_{k+1}$.

\vspace{0.2cm} \noindent
(2) \raisebox{-0.55cm}{\picB{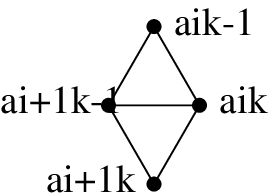}} \hspace{0.1cm} This says
that $\lambda^{(i)}_k \geq \lambda^{(i+1)}_k$. \vspace{0.1cm}

Together, (1) and (2) say that $\lambda^{(i)}_k \geq \lambda^{(i+1)}_k
\geq \lambda^{(i)}_{k+1}$.  In particular, each sequence
$\lambda^{(i)}$ is weakly decreasing, and we have a nested sequence of
partitions: $\lambda^{(1)} \supset
\lambda^{(2)} \supset \dots \supset \lambda^{(n)} \supset
\lambda^{(n+1)} = \emptyset$.

For example, the hive
\[ \picA{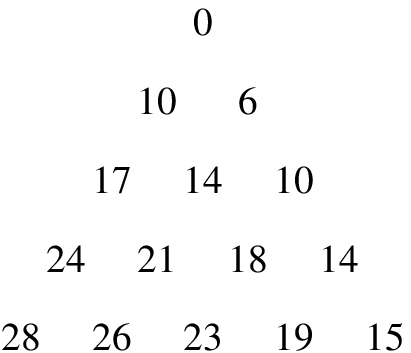} \]
gives the chain of partitions $(6,4,4,1) \supset (4,4,1) \supset (4,2)
\supset (2)$.

We identify partitions with Young diagrams, but rotated by 180
degrees, so the diagram for a partition $\lambda$ has $\lambda_k$
boxes in the $k^\tth$ row from the bottom, and the rows are
lined up on the right.  Fill the boxes by putting the integer $i$ in
each box of $\lambda^{(i)} - \lambda^{(i+1)}$.  The conditions (1) and
(2) say exactly that the result $T$ is a skew tableau on this shape,
that is, it is weakly increasing across rows and strictly increasing
down columns.  Such a $T$ is often called a contratableau of shape
$\lambda$. In our example, $T$ is
\[ \picA{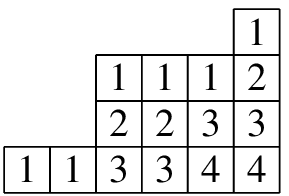} \]

The word $w(T)$ is obtained by reading from left to right in rows, from 
bottom to top.  In the example, $w(T) = 1\, 1\, 3\, 3\, 4\, 4\, 2\, 2\, 3
\,3 \,1 \,1 \,1 \,2 \,1$.

Let $U(\mu)$ be the tableau of shape $\mu$ whose $i^\tth$ row
has $\mu_i$ entries, all equal to $i$.  The word $w(U(\mu))$ is
similarly read from left to right, bottom to top.  In our example,
$\mu = (4,4,3,2)$, and $w(U(\mu)) =
4\,4\,3\,3\,3\,2\,2\,2\,2\,1\,1\,1\,1$.

Now we consider the last rhombus inequalities:

\vspace{0.1cm} \noindent
(3) \raisebox{-0.3cm}{\picB{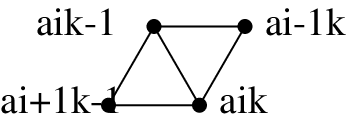}} \vspace{0.15cm}  These
say that $a^{i+1}_{k-1} - a^i_{k-1} \leq a^i_k - a^{i-1}_k$.  We claim
that this is equivalent to the condition that $w(T) \cdot w(U(\mu))$
is a reverse lattice word.

This asserts that, if we divide this word at any point, the number of
times that $i$ occurs to the right of this point is no larger than the
number of times that $i-1$ occurs.  We only need to check this at a
division corresponding to the place in the $k^\tth$ row from the
bottom of $T$ that divides elements strictly smaller than $i$ from
elements greater than or equal to $i$.  The number of times that $i$
occurs here is
\begin{align*}
& (\lambda^{(i)}_k - \lambda^{(i+1)}_k) + (\lambda^{(i)}_{k+1} -
\lambda^{(i+1)}_{k+1}) + \dots + (\lambda^{(i)}_{n+1-i} - 0) + \mu_i \\
& = (\lambda^{(i)}_k + \lambda^{(i)}_{k+1} + \dots + \lambda^{(i)}_{n+1-i})
    - (\lambda^{(i+1)}_k + \lambda^{(i+1)}_{k+1} + \dots +
    \lambda^{(i+1)}_{n-i}) + \mu_i \\
& = (a^i_{n+1-i} - a^i_{k-1}) - (a^{i+1}_{n-i} - a^{i+1}_{k-1}) +
    (a^{i+1}_{n-i} - a^i_{n+1-i}) \\
& = a^{i+1}_{k-1} - a^{i}_{k-1} \,.
\end{align*}
Similarly, the number of times that $i-1$ occurs is 
\[ (\lambda^{(i-1)}_{k+1} - \lambda^{(i)}_{k+1}) + 
   (\lambda^{(i-1)}_{k+2} - \lambda^{(i)}_{k+2}) + \dots + 
   (\lambda^{(i-1)}_{n+2-i} - 0) + \mu_{i-1} = a^i_{k} - a^{i-1}_{k} \,.
\]
Note that the number of times $i$ occurs in all of $T$ is $a^{i+1}_0 -
a^i_0 - \mu_i = \nu_i - \mu_i$.

This process is reversible.  Given any contratableau $T$ of shape
$\lambda$ such that $w(T)\cdot w(U(\mu))$ is a reverse lattice word,
$T$ determines the chain $\lambda^{(1)} \supset \lambda^{(2)} \supset
\dots \supset \lambda^{(n)}$, and from these partitions one
successively fills in the entries in the northwest to southeast
diagonal rows of the hive; the rhombus inequalities (1)--(3) are
automatically satisfied.

To make the story complete, we recall why such contratableaux
correspond to Littlewood-Richardson skew tableaux, using standard
results about tableaux, as in \cite{fulton:young}. However, it may be
pointed out that these contratableaux are at least as easy to produce
and enumerate as the more classical skew tableaux.  First, the
condition that $w(T)\cdot w(U(\mu))$ is a reverse lattice word, given
that the number of times $i$ occurs in $T$ is $\nu_i - \mu_i$, is
equivalent to asserting that $w(T)\cdot w(U(\mu))$ is Knuth equivalent
to $w(U(\nu))$ \cite[\S 5.2]{fulton:young}. The rectification $R$ of a
contratableau $T$ of shape $\lambda$ is easily seen to be a tableau of
shape $\lambda$, and with the same property that $w(R) \cdot
w(U(\mu))$ is Knuth equivalent to $w(U(\nu))$.  The correspondence
between tableaux and contratableaux of shape $\lambda$ is a bijection,
by reversing the rectification process.

Now the condition that $w(R) \cdot w(U(\mu))$ be Knuth equivalent to
$w(U(\nu))$ is equivalent to the condition that $R \cdot U(\mu) =
U(\nu)$ in the plactic monoid of tableaux \cite[\S 2.1]{fulton:young}.
This is easy to see, from the definition of multiplying tableaux by
column bumping entries of the first tableau into the second \cite[\S
A.2]{fulton:young}, that if $R$ and $S$ are tableaux with $R \cdot S =
U(\beta)$, then $S$ must be equal to $U(\alpha)$ for some partition
$\alpha$.  This gives a correspondence between the set of tableaux $R$
that we are looking at and the set of pairs $(R,S)$ with $R$ of shape
$\lambda$, $S$ of shape $\mu$, whose product is the tableau $U(\nu)$.
There is a standard construction \cite[\S 5.1]{fulton:young} between
these pairs and the set of skew tableau on the shape $\nu/\lambda$ of
content $\mu$ whose word is a reverse-lattice word.

%% file: bibliography.tex
\providecommand{\bysame}{\leavevmode\hbox to3em{\hrulefill}\thinspace}

%% file: sat.bbl
\begin{thebibliography}{10}

\bibitem{berenstein.zelevinsky:triple}
A.~D. Berenstein and A.~V. Zelevinsky, \emph{Triple multiplicities for ${\rm
  sl}(r+1)$ and the spectrum of the exterior algebra of the adjoint
  representation}, J. Alg. Comb. \textbf{1} (1992), 7--22.

\bibitem{buch.fulton:chern}
A.~S. Buch and W.~Fulton, \emph{Chern class formulas for quiver varieties}, to
  appear in Invent. Math., 1998.

\bibitem{carre:rule}
C.~Carr{\'e}, \emph{The rule of {L}ittlewood-{R}ichardson in a construction of
  {B}erenstein-{Z}elevinsky}, Internat. J. Algebra Comput. \textbf{1} (1991),
  473--491.

\bibitem{fulton:young}
W.~Fulton, \emph{Young tableaux}, London Mathematical Society Student Texts,
  vol.~35, Cambridge University Press, 1997.

\bibitem{fulton:eigenvalues}
\bysame, \emph{Eigenvalues of sums of {H}ermitian matrices (after {A}.
  {K}lyachko)}, Bourbaki talk, 1998.

\bibitem{horn:eigenvalues}
A.~Horn, \emph{Eigenvalues of sums of {H}ermitian matrices}, Pacific J. Math.
  \textbf{12} (1962), 225--241.

\bibitem{klyachko:stable}
A.~A. Klyachko, \emph{Stable bundles, representation theory and {H}ermitian
  operators}, Institute Mittag-Leffler Preprint, 1996--1997.

\bibitem{knutson.tao:honeycomb}
A.~Knutson and T.~Tao, \emph{The honeycomb model of the
  {B}erenstein-{Z}elevinsky polytope {I}: {K}lyachko's saturation conjecture},
  preprint, 1998.

\bibitem{sturmfels.thomas:variation}
B.~Sturmfels and R.~R. Thomas, \emph{Variation of cost functions in integer
  programming}, Math. Programming \textbf{77} (1997), no.~3, Ser. A, 357--387.

\bibitem{zelevinsky:littlewood-richardson}
A.~Zelevinsky, \emph{Littlewood-{R}ichardson semigroups}, MSRI Preprint,
  1997-044.

\end{thebibliography}
